
\input amstex
\documentstyle{amsppt}

\nologo
\magnification=\magstep1
\hsize=145mm
\vsize=220mm
\hcorrection{-7mm}
\vcorrection{-10mm}
\baselineskip=18pt
 \abovedisplayskip=4pt
 \belowdisplayskip=4pt
\parskip=3pt
\parindent=8mm


\define\almb{\allowmathbreak }
\define\dnl{\newline\newline}
\define\lra{{\longrightarrow}}

\define\nl{\newline}
\define\nomb{\nomathbreak }

\define\Bs{{\roman{Bs}}}

\define\Inf{{\roman{Inf}}}

\define\Ker{{\roman{Ker}}}

\define\Pic{{\roman{Pic}}}
\define\Sing{{\roman{Sing}}}

\define\dg{{\roman{deg}}}
\define\dm{{\roman{dim}}}

\define\SA{\Cal A}

\define\SE{\Cal E}

\define\SI{\Cal I}

\define\SO{\Cal O}

\define\BC{{\Bbb C}}
\define\BP{{\Bbb P}}
\define\BQ{{\Bbb Q}}

\define\BZ{{\Bbb Z}}

\define\Prf{{\it Proof. }}

\define\Rmk{{\it Remark. }}
\define\Th{{\bf Theorem. }}

\define\ke{{\kappa\epsilon}}

\topmatter
\author Takao FUJITA \endauthor
\address 
Takao FUJITA\newline
Department of Mathematics\newline
Tokyo Institute of Technology\newline
Oh-okayama, Meguro, Tokyo\newline
152 Japan
\endaddress
\email
fujita\@math.titech.ac.jp
\endemail
\title On Kodaira energy and adjoint reduction of polarized threefolds
\endtitle

\endtopmatter

\document

\noindent{\bf Introduction}

In this paper, as a continuation of [{\bf F4}], we study and classify polarized threefolds whose Kodaira energies are in the range $\ke<-1/2$.

The Kodaira energy of a polarized manifold $(M, L)$, a pair consisting of a smooth complex algebraic variety $M$ and an ample line bundle $L$ on it, is defined by
$$\ke(M,L):=-\Inf\{t\in\BQ\vert \kappa(K+tL)\ge 0\},$$
where $K$ is the canonical bundle of $M$ and $\kappa$ denotes the Iitaka dimension of a $\BQ$-bundle.
Using the theory of first and second reduction (cf. [{\bf BFS}], [{\bf BS}]), we gave a classification of polarized manifolds with $n=\dm M\ge 4$ and $\ke(M,L)<3-n$ in [{\bf F4}], supplementing the results in [{\bf BS}].
On the other hand, in case $n=3$, we have proved the Spectrum Conjecture (cf. [{\bf F6}]), which says that, for any $\epsilon>0$, there are only finitely many possible values for $\ke(M,L)$ of smooth polarized manifolds $(M,L)$ in the range $\ke<-\epsilon$.
Thus, philosophically, we should be able to classify smooth polarized threefolds with $\ke(M,L)<-\epsilon$ for any $\epsilon>0$.
The purpose of this paper is to carry out such a classification in the range $\ke(M,L)<-1/2$.
The result is summarized at the end (3.$\infty$) of \S3.
In fact, our method here is just a variant of [{\bf F4}].
In order to proceed beyond the range $\ke<-1/2$, we must study more precisely the structure of the birational transformation to obtain canonical fibrations (see [{\bf F6}] for details), which will be left to future investigations.
\dnl
{\bf \S1. Preliminaries}

The results in [{\bf F4};\S1] are our main tools here too.
In the following sections, [{\bf F4};(1.$\ast$)] will be quoted just as (1.$\ast$).
\dnl
{\bf \S2. Second reduction}

In this section also we follow the same line as in [{\bf F4}], except the problem discussed in (2.4).

(2.1) Let $(M,L)$ be a smooth polarized threefold.
We have a classification theory in case $\ke(M,L)\le -1$, so we assume $\ke(M,L)>-1$ from now on.
Moreover, replacing $(M,L)$ by its first reduction if necessary, we assume that $K+L$ is nef.
By the Fibration Theorem (1.2), there are a normal polarized variety $(M'',A)$ and a morphism $f: M\to M''$ such that $f^*A=K+L$ and $f_*\SO_M=\SO_{M''}$.
The assumption $\ke>-1$ implies that $f$ is birational.
This pair $(M'',A)$ will be called the {\it second reduction} of $(M, L)$.
As we shall see later, $M''$ is not always smooth; it may have certain 2-factorial non-Gorenstein terminal singularities.
$f_*L=L''$ is a reflexive sheaf on $M''$, but not always invertible.
Moreover, $L''$ is not always nef as a $\BQ$-bundle.

(2.2) Let $R$ be an extremal ray such that $(K+L)R=0$ and let $\rho: M \to W$ be the contraction morphism of $R$.
As in [{\bf F4};(2.2)], $f$ factors through $\rho$, so $\rho$ is birational.
Let $E$ be the exceptional set of $\rho$.

(2.3) The type of $\rho$ is classified as follows.\newline
1) $\rho(E)=C\subset W$ is a smooth curve, $\rho$ is the blow-up along $C$ and $L=\rho^*L^\flat -E$ for some line bundle $L^\flat$ on $W$, where $E$ is the exceptional divisor of $\rho$.\newline
2) $\rho(E)$ is a point and $(E, L_E)$ is a (possibly singular) hyperquadric in $\BP^3$. Moreover $[E]_E=\SO(-1)$.\newline
3) $\rho(E)$ is a point, $(E, L_E)\simeq(\BP^2, \SO(1))$ and $[E]_E=\SO(-2)$.\newline
4) $\rho(E)$ is a point and $(E, L_E)\simeq(\BP^2, \SO(2))$, $[E]_E=\SO(-1)$.

{\it Remark.\ }
$W$ is smooth in case 1) and 4), but $L^\flat$ is not always nef since $L^\flat C$ may be negative.
In case 2), $\rho(E)$ is a non-factorial hypersurface terminal singularity of $W$, but every Weil divisor on $W$ is Cartier.
In case 3), $\rho(E)$ is a 2-factorial terminal non-Gorenstein singularity and $\rho_*L$ is not invertible.

(2.4) Let $R_1, \cdots, R_k$ be the extremal rays with $(K+L)R=0$ and let $E_j$ be the exceptional set of the contraction morphism $\rho_j: M\to W_j$ of $R_j$.
Unlike the higher dimensional cases (cf. [{\bf F4}; (2.4)]), they may meet with each other.
But we have the following

\Th {\sl Suppose that $E_i\cap E_j\ne\emptyset$ for some $i\ne j$.
Then either $E_i$ or $E_j$ is of the type (2.3;1) with $L^\flat$ ample.}

\Prf
Assume that $\rho_i(E_i)$ and $\rho_j(E_j)$ are points.
$Z=E_i\cap E_j$ is a 1-cycle in $E_i$, and $E_i$ is $\BP^2$ or a hyperquadric, so $Z$ is nef as a divisor on $E_i$, which implies $E_j Z\ge 0$.
On the other hand, $Z$ is $\rho_j$ exceptional, so numerically proportional to $R_j$, hence $E_j Z<0$.
Thus this case is ruled out.

Next suppose that $\rho_i$ is of the type (2.3.1) and $\rho_j(E_j)$ is a point.
We have $L=\rho_i^*L_i-E_i$ for $L_i\in \Pic(W_i)$, and we want to show that $L_i$ is ample on $W_i$.
For this purpose it suffices to show $L_i C_i>0$.
As above, $[E_i]_{E_j}$ is nef and $E_i Z\ge 0$ for $Z=E_i\cap E_j$.
Hence $0<(L+E_i)Z=\rho_i^*L_i\cdot Z=\delta L_i C_i$, where $\delta$ is the degree of the map $Z\to C_i$.
This implies $L_i C_i>0$, as desired.

Finally suppose that both $\rho_i$ and $\rho_j$ are of the type (2.3;1).
If the 1-cycle $E_i\cap E_j$ has a component $Z$ with $E_i Z\ge 0$, then $L_i$ is ample by the above argument.
Similarly, if $E_j Z\ge 0$, then $L_j$ is ample on $W_j$.
Therefore, we may assume that every component of $[E_j]_{E_i}$ has negative self intersection number.
But a $\BP^1$-bundle has at most one curve with negative self intersection, the so-called negative section.
Hence $E_i\cap E_j$ is irreducible, which will be denoted by $Z$.
Note that the line bundle $[-E_i]$ is ample on $E_i$, since $-E_i Z>0$ for the negative section $Z$ and $-E_i F_i>0$ for any fiber $F_i$ of $E_i\to C_i$.

Now we claim $\dm f(E_i)>0$, where $f: M\to M''$ is the second reduction map.
Indeed, if not, $K+L=0$ in $\Pic(E_i)$, so the canonical bundle $K_i$ of $E_i$ is negative since it is the restriction of $-L+E_i$.
This is possible only when $E_i=\Sigma_1$ and $Z$ is the $(-1)$-section on it, but then $-1=K_i Z=-LZ+E_i Z\le -2$, absurd.
Thus the claim is proved.

From this claim we get $f(E_i)=f(Z)$.
Similarly we have $f(E_j)=f(Z)$ by symmetry.
Take a general hyperplane $H$ on $M''$ and let $S=f^*H$.
Then $S$ is a smooth surface, $f_S: S\to H$ is a birational morphism, and there are two $(-1)$-curves over a point on $f(Z)\cap H$ meeting with each other,
which are fibers of $E_i$ and $E_j$.
This is impossible.
Thus the theorem is proved.

(2.5) By virtue of this theorem, the structure of $f: M\to M''$ can be described as follows:
Let \{$R_i$\} and \{$E_i$\} be as above.
If $E_i\cap E_j\ne\emptyset$, we choose $E_i$ as in (2.4), and blow down it first.
Replacing $(M,L)$ by $(W_i,L_i)$, we can continue the same process, since $W_i$ is smooth and $L_i$ is ample.
Thus, finally, we reach a situation where the exceptional sets \{$E_i$\} are disjoint with each other.
From there the second reduction map is just the simultaneous blow down of exceptional sets as in the higher dimensional cases.
Thus, in order to study the structure of $(M'',A)$, it is harmless to assume that $E_i$'s are disjoint with each other from the beginning, hence we pretend so from now on.

In particular we have:\newline
1) {\sl $M''$ has only isolated terminal singularities.}\newline
2) {\sl Any Weil divisor $D$ on $M''$ is Cartier except at $Y_3$, where $Y_3$ is the image of exceptional divisors of the type (2.3;3).
$2D$ is Cartier everywhere on $M''$.}\newline
3) {\sl $L''=f_*L$ is invertible exactly on $M''-Y_3$.
$2L''\cdot Z>0$ for any curve $Z$ such that $Z\not\subset Z_1$, where $Z_1$ is the image of exceptional sets of the type (2.3;1).}

(2.6) As in [{\bf F4}], we decompose $f$ as $M\to M^\sharp \to M''$, where $f^{\sharp}: M\to M^\sharp$ consists of contractions of the types (2.3;1), (2.3;2) and (2.3;4), and $\pi: M^{\sharp} \lra M''$ consists of contractions of the type (2.3.3).
$M^{\sharp}$ will be called the {\it factorial stage}, since every divisor on $M^{\sharp}$ is Cartier.
In particular $L^{\sharp}=f^{\sharp}_*L$ is invertible.

Finally we recall the following

(2.7) {\it Formula. } {\sl Set $\kappa''=\ke(M'', A)$. Then $\ke(M, L)=\kappa''/(1-\kappa'')$.}

For a proof, see [{\bf F4}; (2.7)].
\dnl
{\bf \S3. The structures of second reduction.}

By the general theory in [{\bf F6}], the possible values of Kodaira energies of smooth polarized threefolds has no limit point in the range $\ke<0$.
Moreover, in the range $\ke\le -1$, we have a precise classification theory.
Here we will proceed to the range $-1<\ke<-\frac12$, decide which values are actually occur for $\ke$, and classify and describe the structure of the second reduction.
Our method is a variant of that in [{\bf F4}].

(3.1) Let things be as in \S2, in particular, let $(M'',A)$ be the second reduction of $(M,L)$.
Let $\tau''$ be the smallest number such that $K''+\tau''A$ is nef, where $K''$ is the canonical $\BQ$-bundle of $M''$.
As in [{\bf F4}], by the theory in [{\bf BS}], we have $\tau''=u/2v$ for some positive mutually coprime integers $u, v$ with $u\le 2(n+1)=8$.
By (2.7), we have $\tau''>1$ in case $\ke(M,L)<-\frac12$.
Hence we should consider the following values:
$$\tau''=4, \frac72, 3, \frac52, 2, \frac74, \frac32, \frac43, \frac54, \frac76.$$
In each case we will study the contraction morphism of an extremal ray $R$ such that $(K''+\tau''A)R=0$.

(3.2) The case $\tau''=4$.

$(K''+4A)R=0$ implies $(M'', A)\simeq(\BP^3, \SO(1))$ by (1.6).
So $L''=\SO(5)$ and $\ke=-\frac45$.
This case occurs actually.

(3.3) The case $\tau''=7/2$.

This case is ruled out immediately by (1.6).

(3.4) The case $\tau''=3$.

Again by (1.6), $(M'', A)$ is either a hyperquadric or a scroll over a curve.
Both cases actually occur and $\ke(M, L)=-3/4$.

(3.5) The case $\tau''=5/2$.

We have $(2K''+5A)R=0$, so $2L''R=7AR$.
Let $\rho: M'' \lra W$ be the contraction morphism of $R$.
We divide the cases according to $\dm W$.

(3.5.0) $\dm W=0$. 
We will show that $(M'', A)$ is the projective cone over the Veronese surface $(\BP^2, \SO(2))$.
The proof consists of several steps.

(a) We have $H^q(M'', \SO(tA))=0$ for any $q>0$ and $t\ge -2$ by the vanishing theorem.
So $\chi(M'', tA)=0$ for $t=-1, -2$ and $\chi(\SO_{M''})=1$.
Hence $\chi(M'', tA)=(t+1)(t+2)(dt+3)/6$ for $d=A^3$ by the Riemann-Roch Theorem.

(b) Let $\pi: M^{\sharp} \to M''$ be the factorial stage as in (2.6) and let $E$ be the exceptional divisor of $\pi$, which is the union of components of the type (2.3;3).
Note that $L^{\sharp}=\pi^*L''-\frac12 E$ is Cartier.
Set $B=L^{\sharp}-3A\in\Pic(M^{\sharp})$, where (and also in the sequel) $A$ denotes $\pi^*A$ by abuse of notation.
Then we have $2B=A-E$, so
$\chi(M^{\sharp}, 2tB)=\chi(M^{\sharp}, tA)-\sum_{j=0}^{t-1}\chi(E, [tA-jE]_E)
=\chi(M'',tA)-\mu\sum (2j+1)(2j+2)/2
=(t+1)(t+2)(dt+3)/6-\mu t(t+1)(4t-1)/6$,
where $\mu$ is the number of components of $E$.
This formula is true not only for all integers $t$, but also for all half-integers.
In particular $\chi(M^{\sharp},-B)=(6-d-2\mu)/16$ and
$\chi(M^{\sharp},-3B)=(d+14\mu-2)/16$.

(c) On the other hand, $-B-K^{\sharp}=2A$ for the canonical bundle $K^{\sharp}$ of $M^{\sharp}$, so $H^q(M^{\sharp},-B)\almb =0$ for $q>0$ by the Vanishing Theorem (1.1).
Hence $\chi(M^{\sharp}, -B)=0$ and $d+2\mu=6$.

Moreover, we have $-3B-E-K^{\sharp}=A$, so $H^q(M^{\sharp},-3B-E)=0$ for any $q>0$.
Using the exact sequence
$0\to \SO_{M^{\sharp}}(-3B-E) \to \SO_{M^{\sharp}}(-3B) \to \SO_E(-3) \to 0$,
we infer $\chi(M^{\sharp},-3B)=\mu\chi(\BP^2,\SO(-3))=\mu$,
which yields $2=d-2\mu$ by (b).

From these equalities we obtain $d=4$ and $\mu=1$,
so $\chi(M^{\sharp}, sB)=(s+1)(s+2)/2$ for any $s$.

\comment
(d) We claim that $L^{\sharp}$ is ample on $M^{\sharp}$.

To see this, it suffices to show $L^{\sharp}Z>0$ for any component $Z$ of the set $Z_1$ as in (2.5;3).
This is clear since $Z_1\cap E=\emptyset$ by (2.4) and so $2L^{\sharp}Z=2L''Z=7AZ$.
\endcomment

(d) Since $B+E-K^{\sharp}=3A$, we have $H^q(M^{\sharp},B+E)=0$ for $q>0$.
Using the exact sequence $0 \to \SO_{M^{\sharp}}(B) \to \SO_{M^{\sharp}}(B+E) \to \SO_E(-1) \to 0$,
we get $H^q(M^{\sharp}, B)=0$ for $q>0$ and
$h^0(M^{\sharp}, B)=3$, since $\chi(M^{\sharp}, B)=3$ by the formula in (c).

(e) Now we study the rational map defined by the linear system $\vert B\vert$.
Let $g: \tilde M \to M^{\sharp}$ be a resolution of points of indeterminacy such that $g^*\vert B\vert =D+\Lambda$, where $D$ is the fixed part and $\Lambda$ is the moving part such that $\Bs\Lambda=\emptyset$.
Here, we blow-up only centers mapped into $\Bs\vert B\vert$.
Thus, singular points of $M^{\sharp}$ not in $\Bs\vert B\vert$ survive on $\tilde M$.
Let $\beta: \tilde M \to \BP^2$ be the morphism defined by $\Lambda$, let $Y=\beta(\tilde M)$ be the image and let $X$ be a general fiber of $\tilde M \to Y$.
Set $k=\dm X=3-\dm Y\ge 1$ and $\delta=\dg Y$.
We will compute intersection numbers on $\tilde M$.

Since $A=2B+E=E+2D+2H$ (here $H=\beta^*\SO(1)=[\Lambda]$ and $E$ denotes the total transform $g^*E$), we have $4=A^3\ge 2A^2 H\ge 4AH^2$.

(f) If $k=2$, then $Y$ is a curve of degree $\delta\ge 2$ and $H\sim\delta X$.
Moreover $A^2X=A(2H+E+2D)X=2ADX$ since $AE=0$ and $HX=0$ in the Chow ring of $\tilde M$.
Thus $4=A^3\ge 2A^2H=4\delta ADX$, which implies $ADX=0$ and $A^2X=0$.
This is absurd since nef-bigness is preserved by restriction to general fibers.
Hence we conclude $k=1$, so $\beta$ is surjective.

(g) Now we have $4=A^3\ge2A^2H\ge 4AH^2=4AX>0$, so the equalities hold and $AX=1$.
Take a general member $S$ of $\vert lA\vert$ on $\tilde M$ for $l\gg 0$.
Then $(A-2H)A=(A-2H)^2=0$ on $S$, so $A-2H=E+2D$ is numerically trivial on $S$.
This implies $D\cap S=\emptyset$, hence $\pi(g(D))$ is a finite (possibly empty) subset of $M''$.

(h) We claim that $B$ is nef.
Indeed, if $BZ<0$ for some curve $Z\subset M^{\sharp}$, then $Z\subset\Bs\vert B\vert=g(D)$, so $\pi(Z)$ is a point, hence $Z\subset E$.
But $B_E=\SO(1)$, so this cannot occur.

(i) Thus $A-E-K^{\sharp}=2A+3B$ is nef and big.
Hence $H^1(M^{\sharp},A-E)=0$ and $H^0(M^{\sharp}, A) \to H^0(E, A_E)\cong\BC$ is surjective, so $V\cap E=\emptyset$ for any general member $V$ of $\vert \pi^*A\vert$.
We identify $V$ and $\pi(V)\in\vert A\vert$ from now on.

(j) The restriction $B_V$ to $V$ is ample since $2B=A-E$.
Moreover $B^2V=\nomb1$.
Hence $V$ is irreducible and reduced.

On the other hand, $H^0(M^{\sharp}, B-A)=\nomb0$ since $A^2(B-A)<0$, so we get $h^0(V, B_V)\ge h^0(M^{\sharp}, B)=3$
by using the exact sequence
$0 \to \SO_{M^{\sharp}}(B-A) \to \SO_{M^{\sharp}}(B) \to \SO_V(B) \to 0$.
Hence $(V, B_V)\simeq(\BP^2, \SO(1))$.

Thus $V$ is an ample divisor on the normal variety $M''$, and $(V, [V]_V)\simeq(\BP^2, \SO(2))$.
Hence $M''$ is the projective cone over $V$ by [{\bf B}].

(k) Now we see that $M^{\sharp}$ is the blow-up of $M''$ at the vertex, so $M^{\sharp}\simeq\BP(\SE)$ for the vector bundle $\SE=\SO(2)\oplus\SO$ on $\BP^2$.
$A=\pi^*A$ is the tautological line bundle $H(\SE)=\SO(1)$, and $B$ is the pull-back of $\SO(1)$ on $\BP^2$.
$E$ is the unique member of $\vert A-2B\vert$ and $L^{\sharp}=3A+B$.

\Rmk $(M,L)$ is not always isomorphic to $(M^{\sharp},L^{\sharp})$.

(3.5.1) $\dm W=1$.
In this case let $F$ be a general fiber of $\rho$.
Then $F$ is smooth since $M''$ has only isolated singularities.
Moreover $2K_F+5A_F=0$, but this is ruled out by (1.6).

(3.5.2) $\dm W=2$.
This case is ruled out as above.

(3.5.3) $\dm W=3$.
In this case $\rho$ is birational, but not an isomorphism.
Take a point $x\in W$ with $\dm\rho^{-1}(x)=m>0$.
Then $(K''+(m+1)A)R>\nomb0$ by (1.5), so $m=2$.
Moreover, for the normalization $\tilde X$ of a two-dimensional component $X$ of $\rho^{-1}(x)$, we have $(\tilde X, A_{\tilde X})\cong(\BP^2,\SO(1))$.
Take a general line $\ell$ on $\tilde X$.
Then $A\ell=1$ and $K''\ell\in\BZ$, since $M''$ has only isolated singularities.
But $\ell$ is proportional to $R$, so $(2K''+5A)\ell=0$, contradiction.
Thus this case is ruled out.

(3.6) The case $\tau''=2$.

Let $\rho: M''\to W$ be the contraction morphism of a ray $R$ such that $(K''+2A)R=0$.
We divide the cases according to $\dm W$.

(3.6.0) $\dm W=0$.
Using the Vanishing Theorem we infer that $(M'',A)$ is a Del Pezzo variety.
In particular $K''=-2A$ is invertible and $M''=M^{\sharp}$.
Of course $L''=3A$ and $\ke(M,L)=-\frac23$.
There are several examples of this type.
$M''$ may have some singularities of the type (2.3;2).

(3.6.1) $\dm W=1$.
In this case $\rho$ makes $M''$ a hyperquadric fibration over $W$.

(3.6.2) $\dm W=2$.
In this case $\rho$ makes $M''$ a scroll over $W$.
One easily sees that these cases (3.6.1) and (3.6.2) actually occur.

(3.6.3) $\dm W=3$.
In this case $\rho$ is birational.
We will show that there is a divisor $D$ contained in the smooth locus of $M''$ such that $(D,A_D)\simeq(\BP^2,\SO(1))$ and $[D]_D=\SO(-1)$, and $\rho$ is the contraction of $D$ to a smooth point.
The proof consists of several steps.

(a) Take a point $x\in W$ such that $\dm\rho^{-1}(x)=m>0$.
By (1.5) we infer $m=2$.
Hence $\rho$ is a divisorial contraction.
Let $D$ be the exceptional divisor.
Then, by (1.5), $(\tilde D,A_{\tilde D})\simeq(\BP^2,\SO(1))$ for the normalization $\tilde D$.
In particular $A^2 D=1$.

(b) Let $\pi: M^{\sharp}\to M''$ be the factorial stage as in (2.6) and let $E$ be the exceptional divisor of $\pi$.
$E$ is the union of divisors of the type (2.3;3).
Let $D^\sharp $ be the proper transform of $D$ on $M^\sharp $ and let $E=E'+E''$, where $E'$ is the union of components $E_i$ such that $[D^\sharp ]_{E_i}=\SO(\delta_i)$ with $\delta_i$ being odd,
while $E''$ is the union of $E_j$'s such that $[D^\sharp ]_{E_j}=\SO(\delta_j)$ with $\delta_j$ being even.
Set further $D^+=D^\sharp +\sum'(\frac{\delta_i+1}2)E_i+\sum''\frac{\delta_j}2 E_j$,
where $\sum'$ (resp. $\sum''$) is the sum of components of $E'$ (resp. $E''$).
Note that $[D^+]_{E_i}=\SO(-1)$ for any component $E_i$ of $E'$ and
$[D^+]_{E_j}=\SO$ for any component $E_j$ of $E''$.

(c) We claim $\pi_*\SO(-D^+)=\SI_D:=\Ker(\SO_{M''}\to \SO_D)$.
Indeed, obviously $\pi_*\SO(-D^+)\almb\subset\SI_D$.
On the other hand, since $\SI_D\subset\SO_{M''}=\pi_*\SO_{M^\sharp }$, a local section of $\SI_D$ on $U\subset M''$ corresponds to a function on $\pi^{-1}(U)$ vanishing on $D^\sharp $.
Let $\mu_k$ be the order of its zero along $E_k$ and let $Z$ be the remaining zero locus other than $D^\sharp $ and $E_k$'s.
Then $[D^\sharp +\mu_k E_k+Z]_{E_k}=0$ and $\delta_k-2\mu_k\le 0$ for each $k$, so the above function belongs to $\SO(-D^+)$.
Thus we prove the claim.

(d) We have $tA+E-(K^\sharp +\frac12E)=tA-K''$ in $\Pic(M^\sharp )$ (here and in the sequel we often omit the symbol $\pi^*$) and
$(tA-K'')R\ge0$ for any $t\ge -2$.
Hence $H+tA-K''$ is nef and big for a sufficiently ample line bundle $H$ on $W$.
Therefore $H^q(M^\sharp ,H+tA+E)=0$ for any $t\ge-2$, $q>0$.
Using the exact sequence
$0\to\SO_{M^\sharp }(H+tA) \to \SO_{M^\sharp }(H+tA+E)\to \SO_E(-2)\to 0$,
we get $H^q(M^\sharp ,H+tA)\simeq H^q(M^\sharp ,H+tA+E)$ for any $q$, $t$.

Similarly, we have $H^q(M^\sharp ,H+tA-D^+)\simeq H^q(M^\sharp ,H+tA+E-D^+)$, since
$0\to \SO_{M^\sharp }(H+tA-D^+) \to \SO_{M^\sharp }(H+tA+E-D^+) \to \SO_{E'}(-1)\oplus\SO_{E''}(-2)\to 0$ is exact.
Moreover $(tA-K''-D)R\ge0$ for any $t\ge-2$, so $H+tA+E-D^+ -(K^\sharp +\frac12 E'')$ is nef and big, hence $H^q(M^\sharp ,H+tA+E-D^+)=\nomb0$ for $t\ge-2$, $q>0$.

Now, using
$0 \to \SO_{M^\sharp }(H+tA-D^+) \to \SO_{M^\sharp }(H+tA) \to \SO_{D^+}(tA) \to 0$,
we obtain $H^q(D^+,tA)=\nomb0$ for $t\ge-2$, $q>0$.
This implies $\chi(D^+,tA)=0$ for $t=-2{\text{ and }}-1$, while $\chi(D^+,\SO)=1$.
Hence $\chi(D^+,tA)=(t+2)(t+1)/2$.

(e) Since $[-D^+]_{E_i}=\SO(1)$ for $E_i$ in $E'$ and $[-D^+]_{E_j}=\SO$ for $E_j$ in $E''$,
we have $R^1\pi_*\SO_{M^\sharp }(-D^+)=0$.
Using (c), we infer that
$0 \to \SI_D \to \pi_*\SO_{M^\sharp }(=\SO_{M''}) \to \pi_*\SO_{D^+} \to 0$
is exact, hence $\SO_D=\pi_*\SO_{D^+}$.
Therefore $h^0(D,A)=h^0(D^+,A)=\chi(D^+,A)=3$ by (d), so the $\Delta$-genus of $(D,A_D)$ is zero and $(D,A_D)\simeq(\BP^2,\SO(1))$.

(f) Assume that $p_k:=\pi(E_k)$ is in $D$ for some component $E_k$ of $E$.
A neighborhood $U$ of $p_k$ in $M''$ can be embedded in some $\BC^N$.
Let $\tilde{\BC^N}\to \BC^N$ be the blow-up at $p_k$ and let $G$ be the exceptional divisor $\simeq\BP^{N-1}$.
For the proper transform $\tilde U$ of $U$, we identify $\tilde U\cap G$ with $E_k$, which is a Veronese surface.
On the other hand $\tilde D\cap G$ is a line for the proper transform $\tilde D$ of $D$.
But a Veronese surface contains no line in $G$.
From this contradiction we infer that $D$ does not meet $\pi(E)$.

(g) Now we have $D^+\simeq D$ and $D$ is a Cartier divisor on $M''$.
Therefore $M''$ is smooth along $D$, since $D$ is smooth.
Thus we complete the proof.

{\it Conclusion of (3.6.3).}
$W$ has no worse singularity than $M''$, and $\rho_*A$ is an ample line bundle on $W$.
So, as in [{\bf F4};(4.4.$\infty$)], by the further reduction replacing $(M'',A)$ by $(W,\rho_*A)$, we can get rid of this case.
Note that $\tau''(W,\rho_*A)\le 2$.

(3.7) The case $\tau''=7/4$.

Let $R$ be an extremal ray with $(4K''+7A)R=0$.
Since $L''$ is 2-Cartier, we have $B:=2L''-5A\in\Pic(M'')$.
Moreover $2BR=AR>0$ and $(2K''+7B)R=0$.
This contradicts (1.6).

(3.8) The case $\tau''=3/2$.

Let $\rho: M''\to W$ be the contraction morphism of an extremal ray $R$ with $(2K''+3A)R=0$.
We divide the cases according to $\dm W$.

(3.8.0) $\dm W=0$.

We argue as in [{\bf F4};(4.6.0)].
We have $\chi(M'', -A)=0$ and $\chi(M'', \SO)=\nomb1$ by Vanishing Theorem, so
$\chi(M'', tA)=(t+1)(dt^2+bt+6)/6$ for $d=A^3$ and some $b$.
Calculating the coefficient of $t^2$ by Riemann-Roch Theorem we get $b=\frac54 d$, so $\chi(M'', tA)=(t+1)(4dt^2+5dt+24)/24$.

Let $M^{\sharp}$ be the factorial stage, let $E$ be the exceptional divisor of $\pi: M^{\sharp} \to M''$ and set $B=L^{\sharp}-2A\in\Pic(M^{\sharp})$.
Then $E+2B=A$ and
$\chi(M^{\sharp},2tB)=\chi(M'',tA)-\mu\sum_{s=0}^{t-1}\chi(\BP^2,\SO(2s))
=(t+1)(4dt^2+5dt+24)/24-\mu(4t-1)t(t+1)/6$ as in (3.5.0;b), where $\mu$ is the number of components of $E$.
Applying this formula for $t=-\frac12$, we get
$\chi(M^{\sharp}, -B)=\frac12-\frac{d}{32}-\frac{\mu}{8}$.

On the other hand, we have $\chi(M^{\sharp}, -B)=0$ since $-B-K^{\sharp}=A$ is nef and big on $M^\sharp $.
Thus we get $16=d+4\mu$.
Therefore $\chi(M^{\sharp}, sB)=(s+1)(s+2)((2-\mu)s+3)/6$ and $B^3=2-\mu$.
Moreover $(d,\mu)=(16, 0), (12, 1), (8, 2) {\text{ or }} (4, 3)$.

Before studying each case separately, we compute $h^0(M^\sharp ,B)$.
We have $H^q(M^\sharp ,B+E)=0$ for $q>0$ since $B+E-K^\sharp =2A$ is nef and big.
From the exact sequence
$0 \to \SO_{M^\sharp }(B) \to \SO_{M^\sharp }(B+E) \to \SO_E(-1) \to 0$,
we get $H^q(M^\sharp ,B)=0$ for $q>0$, hence $h^0(M^\sharp ,B)=\chi(M^\sharp ,B)=5-\mu$.

(3.8.0.0) The case $\mu=0$, $d=16$.

In this case $M^{\sharp}=M''$ and $K''=-3B$, so $M''$ is a (possibly singular) hyperquadric in $\BP^4$ and $B=\SO(1)$.
Of course $A=\SO(2)$, $L''=\SO(5)$ and $\ke(M, L)=-3/5$.
There are indeed several such cases.

(3.8.0.1) The case $\mu=1$, $d=12$.

We have $h^0(M^{\sharp}, B)=4$.
Let $g: \tilde M \to M^{\sharp}$ and $g^*{\vert}B{\vert}=D+\Lambda$ be as in (3.5.0.e).
Furthermore, let $X$ be a general fiber of $\beta: \tilde M \to Y\subset\BP^3$ as there and set $k=\dm X=3-\dm Y$ and $\delta=\dg Y$.
For $H=[\Lambda]=\beta^*\SO(1)$ we have $A=E+2B=E+2D+2H$ in $\Pic(\tilde M)$,
so $12=A^3\ge 2A^2H\ge 4AH^2\ge 8H^3$.
Now we divide the cases according to $k$.

(3.8.0.1.0) Suppose that $k=0$.
Then $Y=\BP^3$ and $H^3>0$.
In fact $H^3=1$ by the above inequality.

From the Index Theorem we infer $(A^3/A^2H)^2\le A^2H/H^3$, so $A^2H>\nomb5$.
On the other hand we have $12=A^3=2A^2(D+H)$ since $AE=0$ in the Chow ring.
From them we infer $A^2D=0$ and $A^2H=6$.
Simlarly we have $AH^2>2$ and $6=A^2H=2A(D+H)H$, so $ADH=0$ and $AH^2=3$.
Thus $(E+2D)H^2=AH^2-2H^3=1$, hence $EH^2=1$ and $DH^2=0$.

Now, by the same argument as in [{\bf F4};(4.6.0.1.0;b \& c)], we infer $D=0$.
Thus $\Bs{\vert}B{\vert}=\emptyset$, $\tilde M=M^{\sharp}$, $H=B$, $\beta: M^{\sharp} \to \BP^3$ is a birational morphism and $\beta(E)$ is a hyperplane.
Hence $\beta^*(\beta(E))=E+E^*$ for some member $E^*$ of ${\vert}B-E{\vert}$.
We will show that $\beta$ is the blow-up along $C=\beta(E^*)$.

We have $B^2E^*=B^3-B^2E=0$, so $\dm C\le 1$.
On the other hand $ABE^*=AB^2-ABE=3$, so $C$ is a curve in $\beta(E)\simeq\BP^2$.
We have $H^1(M^{\sharp}, A-E)=0$ since $A-E-K^{\sharp}=A+3B$ is nef and big.
Hence $H^0(M^\sharp ,A)\to H^0(E,A_E)$ is surjective and $E\cap T=\emptyset$ for any general member $T$ of ${\vert}A{\vert}$.
Since $[T]=A=3B-E^*$, $\beta(T)$ is a hypercubic in $\BP^3$ containing $C$.
Therefore $C=\beta(T)\cap\beta(E)$, $\beta^*(\beta(T))=T+E^*$, $\beta^*(\beta(E))=E+E^*$ and $T\cap E=\emptyset$.
By the universal property of the blowing-up, this implies that $\beta$ factors through the blow-up $\tilde P$ of $\BP^3$ along $C$, and $E^*$ is the pull-back of the exceptional divisor $E_C$ lying over $C$.
The morphism $M^{\sharp} \to \tilde P$ is finite since the ample line bundle $L^{\sharp}=2A+B=7B-2E^*$ comes from $\Pic(\tilde P)$.
Moreover, as in [{\bf F4};(4.6.0.1.0;f)], $\tilde P$ has only hypersurface singularities and ${\roman{codim}}\Sing(\tilde P)\ge 2$, so $\tilde P$ is normal.
Therefore $M^{\sharp}\simeq\tilde P$ by Zariski's Main Theorem.

The situation can be described in the following way too.
Let $G$ be the scroll over $\BP^3$ associated with the bundle $\SO(3)\oplus\SO(1)$, and let $H_\gamma$ be the tautological bundle on it, while $H_\beta$ is the pull-back of $\SO_{\BP^3}(1)$.
Let $\Delta_{\infty}$ be the unique member of $\vert H_\gamma-3H_\beta\vert$ and let $\Delta_0\in\vert H_\gamma-H_\beta\vert$ be another section disjoint from $\Delta_{\infty}$.
Then $\tilde P=M^\sharp $ is embedded in $G$ as a member of $\vert H_\gamma\vert$ in such a way that $E=\Delta_{\infty}\cap M^\sharp $ and $T=\Delta_0\cap M^\sharp $.
Moreover $A=H_\gamma-H_\beta$ and $L^\sharp =2H_\gamma-H_\beta$ in $\Pic(M^\sharp )$.
The defining equation of $\tilde P$ in $G$ belongs to $H^0(G,H_\gamma)\simeq H^0(\BP^3,\SO(3)\oplus\SO(1))$,
so $\phi\oplus\psi$ with $\phi\in H^0(\BP^3,\SO(1))$ and $\psi\in H^0(\BP^3,\SO(3))$,
and $C$ is the complete intersection $\{\phi=\psi=0\}$ in $\BP^3$.

Clearly such a case occurs really.
Unlike [{\bf F4};(4.6.0.1.0)], $C$ may have singularities since $L^\sharp Z=2$ for a curve $Z$ in $M^\sharp $ lying over $x\in\BP^3$.

(3.8.0.1.1) Suppose that $k=1$.
Then $Y$ is a surface of degree $\delta\ge2$ in $\BP^3$.
Hence $12=A^3\ge4AH^2=4\delta AX$ implies $1=AX=(E+2D)X$, so $DX=0$ and $EX=1$.
Therefore $\tilde EX=1$ for the proper transform $\tilde E$ of $E$ on $\tilde M$, since every $g$-exeptional component is a component of $D$.
Moreover $[B-H]_{\tilde E}=D_{\tilde E}$ is effective since $\tilde E$ is not a component of $D$.
$B_{\tilde E}$ is nef since $B_E=\SO(1)$, so we infer $H^2\tilde E\le B^2\tilde E=B^2E=1$.
But the restriction of $\beta$ to $\tilde E$ is a birational morphism onto $Y$, so $H^2\tilde E=\delta\ge2$.
Thus this case is ruled out.

(3.8.0.1.2) Suppose that $k=2$.
$Y$ is a curve of degree $\delta\ge 3$ since it is not contained in any plane.
From $12=A^3\ge 2A^2H=2\delta A^2X=4\delta ADX$ we infer $ADX=1$, $\delta=3$ and $A^2X=2$.
Moreover $D$ is mapped onto a curve in $M''$ since $0=A^2(A-2H)=2A^2D$.
We can derive a contradiction by the same argument in [{\bf F4};(4.6.0.1.2)].

(3.8.0.2) The case $\mu=2$, $d=8$.

This time $E$ is the sum of two components $E_1$ and $E_2$.
Let $\tilde M$, $g$, $D$, $H$ be as above.
Here $h^0(M^\sharp ,B)=5-\mu=3$, so we have $\beta: \tilde M \to \BP^2$.
Let $Y=\beta(\tilde M)$, $X$ and $k=\dm X$ be as before.
Clearly $k>0$.

Suppose that $k=2$.
Then $Y$ is a curve of degree $\delta\ge 2$.
Since $8=A^3\ge 2A^2H=2\delta A^2X=2\delta A(E+2D+2H)X=4\delta ADX$,
we infer $ADX=\nomb 1$, $A^2X=2$ and $\delta=2$.
This case is ruled out again by the same method in [{\bf F4};(4.6.0.1.2)].

Now we conclude $k=1$.
We have $8=A^3\ge 2A^2H\ge 4AH^2=4AX\ge 4$.
Since $8=A^3=2A^2B=2A^2(D+H)$ and $A^2D=2ABD$, if $A^2D>0$ we would have $A^2H\le 2$, contradicting the Index Theorem.
Hence $A^2D=0$ and $4=A^2H=2ABH=2A(D+H)H$.

If $ADH>0$, then $ADH=1=AH^2=(E+2D)X$, so $EX=1$ and $DX=0$.
But this implies $D_jX=0$ for any $g$-exceptional divisor $D_j$, so ${\tilde K}X=K^\sharp X$ for the canonical bundle $\tilde K$ of $\tilde M$,
hence $2g(X)-2=\tilde K X=(-\frac32 A+\frac12 E)X=-1$, absurd.
Thus we conclude $ADH=0$ and $AH^2=AX=2$.

Now we claim that $\pi(g(D))$ is at most finite.
To see this, let $S$ be a general member of $\vert \ell A\vert$ with $\ell\gg 0$.
We have $0=A(A-2H)^2=4AD^2$, so $D_S^2=0$ and $AD_S=0$, hence $D_S$ is numerically trivial by the Index Theorem.
This implies $D\cap S=\emptyset$, so $\dm\pi(g(D))\le 0$ since $A$ is ample on $M''$.

As in (3.5.0; h), from this claim we infer that $B$ is nef.
Next we claim $DH^2=0$.
Indeed, otherwise, we have $DH^2=1$ and $EH^2=0$ since $2=AH^2=(E+2D)H^2$.
Let $D_0$ be the unique component of $D$ with $D_0 H^2=1$.
This is not a component of the total transform of $E$, so $g(D_0)$ is a point off $E$, and is an isolated base point of $\vert B\vert$.
Since $B$ is nef, the existence of isolated base point implies $B^3>0$, which contradicts $B^3=2-\mu$.
Thus the claim $DH^2=0$ is proved.

Finally we claim $D=0$.
If not, $\beta(D)$ is a curve by the above claim.
Clearly $D$ is not numerically proportional to $H$, so
$\beta^*\beta(D)=D+D^*$ for some $D^*\ne 0$.
Then we have a curve $Z$ in $D$ such that $\beta(Z)$ is a point and $D^*Z>0$, $DZ<0$.
But this contradicts the nefness of $B=D+H$.

Now we have $\tilde M=M^\sharp $, $\Bs{\vert}B{\vert}=\emptyset$, $B=H$ and the morphism $\beta: M^{\sharp} \to \BP^2$.
For each component $E_i$ of $E$, we have $B_{E_i}=\SO(1)$ and $E_i$ yields a section of $\beta$.
Moreover $E$ is $\beta$-ample since $L^{\sharp}=2A+B=5B+2E$ is ample on $M^{\sharp}$.
Hence every fiber $X$ of $\beta$ is a curve, and has at most two components.
General fiber is $\BP^1$ since $K^{\sharp}X=-2$, so $\beta$ is a conic bundle.

As in [{\bf F4}; (4.6.0.2.1;f)], taking $\beta_*$ of the exact sequence $0 \to \SO_{M^{\sharp}}(2B) \to \SO_{M^{\sharp}}(A) \to \SO_{E_1}\oplus\SO_{E_2} \to 0$,
we get an exact sequence
$0 \to \SO_Y(2) \to \SA \to \SO_Y\oplus\SO_Y \to 0$ for $\SA=\beta_*\SO_{M^{\sharp}}(A)$.
This must split and
$\SA\simeq\SO(2)\oplus\SO\oplus\SO$ on $Y\simeq\BP^2$.
Moreover, as in [{\bf F4}], we get a morphism $\alpha: M^{\sharp} \to P=\BP(\SA)$ such that $\alpha^*\SO_P(1)=A$.
Set $H_{\alpha}=\SO_P(1)$ and let $H_{\beta}$ be the pull-back of $\SO_Y(1)$.
Then $E=\alpha^*T$, where $T$ is the unique member of ${\vert}H_{\alpha}-2H_{\beta}{\vert}$ on $P$ such that $T\simeq Y\times\BP_\alpha^1$.
$\alpha$ is an embedding and $M^{\sharp}$ is a member of ${\vert}2H_{\alpha}{\vert}$ on $P$ since $M^{\sharp}\cap T=E_1+E_2\in{\vert}2H_{\alpha}{\vert}_T$.

Conversely, starting from a general member of ${\vert}2H_{\alpha}{\vert}$ on $P=\BP_{\BP^2}(\SO(2)\oplus\SO\oplus\SO)$, we can construct a polarized threefold of this type.

(3.8.0.3) The case $\mu=3$, $d=4$.
Unlike [{\bf F4}; (4.6.0.3)], this case really occurs.

Note first that $h^0(M^{\sharp},B)=5-\mu=2$ and
$B^3=2-\mu=-1$,
so $B$ is not nef.
Let things be as before.
We have $4=A^3=2A^2B=2(A^2D+A^2X)\ge 2A^2X$,
$A^2X=2ABX=2ADX$ and $A^2X>0$.
Hence $ADX=1$, $A^2X=2$ and $A^2D=0$, which implies $\dm\pi(g(D))\le 1$.
Therefore every fixed component of $\vert B \vert$ is $\pi$-exceptional, so the fixed part of $\vert B\vert$ is of the form $\sum_{i=1}^3 m_i E_i$,
where $E_i$'s are the components of $E$.

Let $B_1+\sum m_iE_i$ and $B_2+\sum m_iE_i$ be general members of $\vert B\vert$.
Then $g(D)=\Bs\vert B\vert=(B_1\cap B_2)\cup(\cup_{m_i>0}E_i)$.
Since $AE_i=0$ in the Chow ring, we have $AB_1B_2=AB^2=\frac14 A(A-E)^2=1$.

Since $A$ is ample on $M''$, $\pi(B_1\cap B_2)$ is an irreducible curve.
Let $C$ be the proper transform of it on $M^\sharp $.
Then $B_1\cdot B_2=C$ modulo 1-cycles contained in $E$, and $AC=1$.

Since $B$ is not nef, there is a curve $Z$ with $BZ<0$.
Such a curve $Z$ must be in $\Bs\vert B\vert$, but $Z\not\subset E$ since $B_{E_i}=\SO(1)$.
Hence $Z=C$.
Thus, $C$ is the unique curve on $M^\sharp $ with $BC<0$.
On the other hand, $L^\sharp $ is ample as in [{\bf F4}; (3.9.0; d)], hence
$0<L^\sharp  C=(2A+B)C=2+BC$, so $BC=-1$.

Now we infer $-1=BC\le BB_1 B_2=B(B-\sum_i m_i E_i)^2=B^3-2B^2\sum m_i E_i+\sum m_i^2BE_i^2=-1-2\sum m_i-2\sum m_i^2$.
This implies $m_i=0$, namely $\vert B\vert$ has no fixed component and $B_1, B_2\in \vert B\vert$.
Moreover $C$ is exactly the scheme-theoretic intersection $B_1\cap B_2$.

Since $K^\sharp =-\frac32 A+\frac12 E=-A-B$, we have $2g(C)-2=(K^\sharp +2B)B^2=-2$.
Hence $C\simeq\BP^1$, so $B_1$ and $B_2$ are smooth along $C$ and intersect normally.
From this we infer that $\tilde M$ is the blow-up of $M^\sharp $ along $C$ and $D$ is the exceptional divisor over $C$.
The normal bundle of $C$ is $B\oplus B\simeq\SO(-1)\oplus\SO(-1)$ and $D\simeq\BP_\xi^1\times C$.
The pull-back of $\SO(1)$ on $C\simeq\BP^1$ will be denoted by $H_\alpha$ from now on.
Thus the conormal bundle $[-D]_D$ is $H_\xi+H_\alpha$.
$B_1\cap E_i$ and $B_2\cap E_i$ are lines on $E_i\simeq\BP^2$ since $B_{E_i}=\SO(1)$.
So $B_1\cap B_2\cap E_i$ is a simple point, hence $C$ and $E_i$ meet normally there.
The proper transform $\tilde{E_i}$ on $\tilde M$ is isomorphic to $\Sigma_1$.

We claim $\Bs\vert A\vert=\emptyset$.
To see this, first note that $\tilde K=K^\sharp +D=-A-B+D=-A-H$, so
$H^q(\tilde M,2H)=0$ for $q>0$ by the Vanishing Theorem.
We have $[D+2H]_D=-H_\alpha+H_\xi$ on $D\simeq\BP_\alpha^1\times\BP_\xi^1$, hence
$H^q(D,D+2H)=0$ for any $q$.
Combining them we get $H^q(\tilde M, D+2H)=0$ for $q>0$
in view of
$0 \to \SO_{\tilde M}(2H) \to \SO_{\tilde M}(D+2H) \to \SO_D(D+2H) \to 0$.
The restriction of $E+D+2H$ to $\tilde{E_i}\simeq\Sigma_1$ is $-D_{\tilde{E_i}}$ and $D\cap\tilde{E_i}$ is the $(-1)$-curve, hence
$H^q(\tilde{E_i},E+D+2H)=0$ for any $q$ and $H^q(\tilde M,E+D+2H)=0$ for $q>0$.
Finally $A=E+2D+2H$ and $A_D=H_\alpha$, so $H^0(\tilde M,A)\to H^0(D,A_D)$ is surjective.
This implies $\Bs\vert A\vert\cap D=\Bs\vert A_D\vert=\emptyset$,
hence $\Bs\vert A\vert\subset E$, because $E+2D+2H\in\vert A\vert$.
For each $i$, by the above surjectivity there is a member $G$ of $\vert A\vert$ such that $G\cap D\cap E_i=\emptyset$.
Since the restriction of $A$ to $E_i$ is trivial, this implies $G\cap E_i=\emptyset$.
Thus $E_i\cap\Bs\vert A\vert=\emptyset$, which implies $\Bs\vert A\vert=\emptyset$ as claimed.

Set $\SA:=\beta_*\SO_{\tilde M}(A)$.
By this claim $\beta^*\SA\to\SO_{\tilde M}(A)$ is surjective, so we have a morphism
$\alpha: \tilde M \to P:=\BP(\SA)$ with $\alpha^*H(\SA)=A$,
where $H(\SA)$ is the tautological line bundle on the scroll $P$.
We next claim $\SA\simeq\SO(2)\oplus\SO\oplus\SO$.

From the exact sequences
$0\to \SO_{\tilde M}[2H] \to \SO_{\tilde M}[D+2H] \to \SO_D[-H_\alpha+H_\xi] \to 0$ and
$0\to \SO_{\tilde M}[D+2H] \to \SO_{\tilde M}[{\tilde E}+D+2H] \to \SO_{\tilde E}[{\tilde E}+D+2H] \to 0$
we infer
$\SO(2)\simeq\beta_*\SO_{\tilde M}(D+2H)\simeq\beta_*\SO_{\tilde M}({\tilde E}+D+2H)$.
Hence from
$0\to \SO_{\tilde M}({\tilde E}+D+2H)\to \SO_{\tilde M}(A) \to \SO_D(A_D) \to 0$
we get an exact sequence
$0 \to \SO(2) \to \SA \to \SO\oplus\SO \to 0$
on $\BP_\xi^1$, so $\SA\simeq\SO(2)\oplus\SO\oplus\SO$, as claimed.

The subbundle $\SO(2)$ of $\SA$ corresponds to the unique member $\Delta$ of $\vert H(\SA)-2H_\xi\vert$ on $P$, where $H_\xi$ is the pull-back of $\SO(1)$ on $\BP_\xi^1$.
Moreover $D\simeq\Delta$ via $\alpha$.
For each $i$, $\alpha(\tilde{E_i})=Y_i$ is a section of $\Delta\to\BP_\xi^1$.
Over them $\alpha$ has fibers of positive dimension, but $\alpha$ is finite over $P-\alpha(\tilde E)$, since every curve $Z$ with $AZ=0$ lies in $\tilde E$ or $D$.
For a general fiber $X$ of $\beta$ we have ${\tilde K}_X=-A_X$ and $A_X^2=2$, so $\alpha_X: X\to\BP^2$ is of degree two.
Hence $\alpha$ is of degree two, and $\alpha$ is ramified over $\Delta$, since $\alpha^*\Delta=E+2D$.

Since $H^0(P,3H(\SA))\to H^0(\Delta,3H(\SA)_\Delta)$ is surjective,
there is a member $T\in\vert 3H(\SA)\vert$ such that $T_\Delta=Y_1+Y_2+Y_3$.
Then $\tilde E_i$'s are components of $\alpha^* T$, so $\alpha^* T=\tilde E_1+\tilde E_2+\tilde E_3+\tilde T$ for some $\tilde T\in\vert 3A-\tilde E\vert$.
We have $\tilde T\cap D=\emptyset$ since $(3A-\tilde E)_D=\SO$.
Thus the scheme theoretic inverse image of $Y_i$ is $\tilde E_i$.
Hence there is a morphism $\tilde\alpha: \tilde M\to \tilde P$ onto the blow-up $\tilde P$ of $P$ along $Y=\cup Y_i$ such that ${\tilde\alpha}^*\Gamma_i=E_i$, where $\Gamma_i$ is the exceptional divisor lying over $Y_i$.

This map $\tilde\alpha$ is a finite morphism of degree two and the branch locus is of the form $\tilde\Delta+R$, where $\tilde\Delta$ is the proper transform of $\Delta$.
The singularities of $\tilde M$ lie exactly over the singularities of $\tilde\Delta+R$, so $\tilde\Delta\cap R=\emptyset$ since $\tilde M$ is smooth along $D={\tilde\alpha}^{-1}(\tilde\Delta)$.
We have $\tilde\Delta+R\in\vert 2F\vert$ for some $F\in\Pic(\tilde P)$, and the canonical bundle $K_{\tilde M}$ of $\tilde M$ is ${\tilde\alpha}^*(K_{\tilde P}+F)$.
Since $K_{\tilde M}=-A-H_\xi$ and $K_{\tilde P}=-3H(\SA)+\sum \Gamma_i$,
we infer $F=2H(\SA)-H_\xi-\sum \Gamma_i$ and $R\in\vert 3H(\SA)-\sum \Gamma_i\vert$.
If $\tilde M$ has a singularity of the type (2.3;2), $R$ has the corresponding singularity at its image on $\tilde P$.

Conversely, starting from a member $R$ of $\vert 3H(\SA)-\sum \Gamma_i\vert$ as above, we can construct a manifold $(M,L)$ of this type.
To see this, we assume that $R$ is smooth for the sake of simplicity.
Let $\tilde\alpha: \tilde M\to \tilde P$ be the double covering branched along $\tilde\Delta+R$, and set
$A={\tilde\alpha}^*H(\SA)$, $2D={\tilde\alpha}^*\tilde\Delta$ and $\tilde E_i={\tilde\alpha}^*\Gamma_i$.
The normal bundle $N_{Y_i\subset P}$ of $Y_i$ is $\SO\oplus\SO(-2)$ since
$0 \to N_{Y_i\subset \Delta} \to N_{Y_i\subset P} \to [N_{\Delta\subset P}]_{Y_i} \to 0$ is exact,
so $\Gamma_i\simeq\Sigma_2$ and $\tilde\Delta\cap\Gamma_i$ is the $(-2)$-section,
while $R\cap\Gamma_i$ is a $(+2)$-section disjoint from it.
Since $\tilde E_i\to\Gamma_i$ is a double covering branched along these sections, we infer $\tilde E_i\simeq\Sigma_1$ and $\tilde E_i\cap D$ is the $(-1)$-curve on it.
Moreover, $[-\tilde E_i]_{\tilde E_i}$ is two times of a $(+1)$-section, since $[-\Gamma_i]_{\Gamma_i}$ coresponds to a $(+2)$-section.
Now, $D\simeq\tilde\Delta\simeq\BP_\xi^1\times\BP_\alpha^1$ and
$[2D]_D={\tilde\alpha}^*[\tilde\Delta]_{\tilde\Delta}=H_\alpha-2H_\xi-3H_\alpha$, so $[D]_D=-H_\alpha-H_\xi$, hence $D$ can be blown down to the direction
$D\to\BP_\alpha^1=C$.
Let $g: \tilde M\to M^\sharp$ be the blow down.
Its restriction to $\tilde E_i$ is the blow down of the $(-1)$-curve $\tilde E_i\cap D$, so $E_i=g(\tilde E_i)\simeq\BP^2$.
Moreover $[E_i]_{E_i}=\SO(-2)$.
We see $A$ and $B=D+H_\xi$ come from $\Pic(M^\sharp)$, so we set
$L^\sharp =2A+B\in\Pic(M^\sharp )$.
By using Nakai's criterion as in [{\bf F4};(2.3)],
we can check that $L^\sharp $ is ample.
Computing the canonical bundles from the side of $\tilde P$, we see that $(M^\sharp ,L^\sharp )$ is a polarized manifold of the desired type.

(3.8.1) The case $\dm W=1$.

Since $M''$ has only isolated singularities,
any general fiber $F$ of $\rho$ is smooth, and the restriction $L''_F$ of $L''$ is invertible.
Set $B=(L''-2A)_F\in\Pic(F)$.
Then $A_F=2B$ and the canonical bundle of $F$ is $K''_F=-\frac32A_F=-3B$.
Hence $(F, B)\simeq(\BP^2, \SO(1))$, $A_F=\SO(2)$ and $L''_F=\SO(5)$.

Next we will show that any singular fiber $X$ of $\rho$ is of the following type:\nl
{\sl $X$ is irreducible, $(X,A)$ is isomorphic to the cone over a rational curve of degree four.
$M''$ has a singularity of the type (2.3;3) at the vertex of $X$, and is smooth elsewhere.
On the factorial stage $M^\sharp$, the proper transform $\tilde X$ is isomorphic to $\Sigma_4$, which can be blown down smoothly to a quadric in $\BP^2$, and the result is a $\BP^2$-bundle near this fiber.}

The proof consists of several steps.
Note that the restriction of the canonical $\BQ$-bundle $K''$ to $X$ is numerically equivalent to $-\frac32 A_X$.
Hence $K''AX_i=-\frac32 A^2X_i\in \BZ$ for any irreducible component $X_i$ of $X$, so $A^2X_i$ is even.
Since $A^2F=4$, $X$ must be of the form $X_1+X_2$ or $2X_1$ if it is not irreducible and reduced.

(a) Suppose that $X=2X_1$.
Let $X_1^{\sharp}$ be the proper transform of $X_1$ on the factorial stage $M^{\sharp}$.
Then $\pi^*X=2X_1^{\sharp}+\sum m_iE_i$ where $E_i$'s are $\pi$-exceptional divisors lying over points on $X$.
The canonical bundle $K_{X_1^{\sharp}}$ of $X_1^{\sharp}$ is $K''+\frac12E+X_1^{\sharp}$ where $E=\sum E_i$, so
$2g(X_1^{\sharp},A)-2=(A+K''+X_1^{\sharp})A\{X_1^{\sharp}\}=(A+K''-\frac12\sum m_iE_i)A\{X_1^{\sharp}\}=(A+K'')A\{X_1^{\sharp}\}=-\frac12 A^2\{X_1^{\sharp}\}=-1$, absurd.
Thus this case is ruled out.

(b) Suppose that $X=X_1+X_2$.
Let $X_i^{\sharp}$ be the proper transform of $X_i$ on $M^{\sharp}$ and let
$\pi^*X=X_1^{\sharp}+X_2^{\sharp}+\sum m_i E_i$ as above.
$X_i^{\sharp}$ is Cartier on $M^{\sharp}$ and its canonical bundle is numerically
$-\frac32A+\frac12 E+X_i^{\sharp}$, so
$2g(X_1^{\sharp},A)-2=(-\frac32 A+\frac12 E+X_1^{\sharp}+A)A\{X_1^{\sharp}\}=(-\frac12 A-X_2^{\sharp})A\{X_1^{\sharp}\}\le-\frac12 A^2\{X_1^{\sharp}\}=-1$,
which yields $g(X_1^{\sharp},A)=0$ and $AX_1^{\sharp}X_2^{\sharp}=1$.

Let $\tilde{X_i^{\sharp}}$ be the normalization of $X_i^{\sharp}$.
Then $0\le g(\tilde{X_1^{\sharp}},A)\le g(X_1^{\sharp},A)$, hence $g(\tilde{X_1^{\sharp}},A)=0$, which yields $\Delta(\tilde{X_1^{\sharp}},A)=0$ by [{\bf F2};(2.11)].
Let $\tilde{X_i}$ be the normalization of $X_i$.
Then $\Delta(\tilde{X_1},A)=\Delta(\tilde{X_1^{\sharp}},A)=0$, hence
$\tilde{X_1}$ is a smooth hyperquadric $\BP^1\times\BP^1$ or a singular hyperquadric, i.e., the projective cone over a plane quadric curve, since $A$ is ample on $\tilde{X_1}$.
In the former case, by lifting a general line on $\tilde{X_1}$, we find a curve $Z$ on $X_1^{\sharp}$ such that $AZ=1$ and $Z\cap E=\emptyset$.
But the restriction of the canonical bundle of $M^{\sharp}$ is $-\frac32 A+\frac12 E$ and we should have $(-\frac32 A+\frac12 E)Z\in\BZ$, contradiction.
Thus $\tilde{X_1}$ is a quadric cone.
Similarly, $\tilde{X_2}$ is also a quadric cone.
Moreover, from $AX_1^{\sharp}X_2^{\sharp}=1$ we infer that $\pi(X_1^{\sharp}\cap X_2^{\sharp})=X_1\cap X_2$ is the image of a line on the cone $\tilde{X_i}$ passing the vertex.

We have $H^0(X,A)\subset\Ker(\psi)$, where $\psi: H^0(X_1,A)\oplus H^0(X_2,A) \to H^0(X_1\cap X_2,A)$ is defined by
$\psi(\phi_1\oplus\phi_2)=-\phi_1\vert_{X_1\cap X_2}+\phi_2\vert_{X_1\cap X_2}$.
Hence we have injections $\Ker(H^0(X,A)\to H^0(X_2,A)) \to \Ker(H^0(X_1,A)\to H^0(X_1\cap X_2,A)) \to \Ker(H^0(X_1,A)\to H^0(L,A))$, where
$L$ is the line on $\tilde{X_1}$ mapped onto $X_1\cap X_2$.
Therefore $h^0(X,A)\le h^0(\tilde{X_2},A)+(h^0(\tilde{X_1},A)-h^0(L,A))=6$.
On the other hand $h^0(X,A)\ge h^0(\BP^2,\SO(2))=6$ by the upper semicontinuity theorem, so the equalities must hold.
In particular $h^0(X_i,A)=h^0(\tilde{X_i},A)$ and $\Delta(X_i,A)=0$, hence $X_i=\tilde{X_i}$.

By an argument as in (3.6.3;f), there is no singularity of the type (2.3;3) on $X_i$ except at the vertex.
Hence $X_i^{\sharp}$ is isomorphic to the blow-up of $X_i$ at the vertex, so $\simeq\Sigma_2$.
Let $Y$ be the fiber contained in $X_1^{\sharp}\cap X_2^{\sharp}$.
Then $X_1^{\sharp}Y\ge 0$ since $Y$ is nef on $X_2^{\sharp}$.
Similarly $X_2^{\sharp}Y\ge 0$.
Moreover $EY>0$, since the restriction of $E$ is the $(-2)$-curve on $\Sigma_2$.
This contradicts $0=\pi^*X\cdot Y=(X_1^{\sharp}+X_2^{\sharp}+mE)Y$.
Thus this case is ruled out.

(c) Now we conclude that $X$ is irreducible and reduced.
Moreover $h^0(X,A)\ge 6$ by the upper semicontinuity theorem, so $\Delta(X,A)=0$, hence $X$ is the projective cone over a Veronese curve of degree four.
The vertex $v$ is the unique singularity of the type (2.3;3) lying on $X$, and the proper transform $X^{\sharp}$ on $M^{\sharp}$ is the blow-up of $X$ at $v$, so $\simeq\Sigma_4$.
$C:=E\cap X^{\sharp}$ is the $(-4)$-curve on it.
Since $EC=-4$, $C$ is a quadric curve in $E\simeq\BP^2$.
We see $\pi^*X=X^{\sharp}+E$, so $[X^{\sharp}]_{X^{\sharp}}=-C$ in $\Pic(X^{\sharp})$, hence $X^{\sharp}$ can be blown down to the direction $X^{\sharp}\to C$.
By this blow down $M^{\sharp}$ is transformed to a $\BP^2$-bundle near this fiber.
Thus the assertion is proved.

(3.8.2) The case $\dm W=2$.

Any general fiber is a smooth curve such that $K_F=-3B$.
This cannot occur.

(3.8.3) The case $\dm W=3$. $\rho$ is birational.

Suppose that $\rho$ is a small contraction.
Then, by [{\bf Ben}], $-K''Z<1$ for a $\rho$-exceptional curve $Z$.
But this is impossible since $(K''+\frac32 A)Z=0$ and $AZ\in\BZ$.
Thus $\rho$ is a divisorial contraction.

Let $D$ be the $\rho$-exceptional divisor.
If $\dm\rho(D)>0$, take a general hyperplane $H$ on $W$ and put $S=\rho^*H$.
Then $S$ is smooth since $\Bs\vert H\vert=\emptyset$ and $M''$ has only isolated singularities.
As a component of $D\cap S$, we find a curve $Z$ such that $\rho(Z)$ is a point and $K_S Z=K''Z=-\frac32 AZ\le -\frac32$.
This yields a contradiction.
Thus $\rho(D)$ is a point.

We proceed as in [{\bf F4};(4.6.4)].
As there, let $D^{\sharp}$ be the proper transform of $D$ on $M^{\sharp}$ and let $D^+$, $D^-$ be as in (3.6.3.c).
Thus $E=E'+E''$, $D^+=D^-+E'$, $\pi^*D=D^+-\frac12E'=D^-+\frac12E'$ and $[D^-]_{E_i}=\SO(1)$ for each component $E_i$ of $E'$, while $[D^-]_{E_j}=\SO$ for each component $E_j$ of $E''$.

We see that $tA-K''$ is $\rho$-ample for $t\ge -1$, hence
$H+E+tA-(K^{\sharp}+\frac12E)=\pi^*(H+tA-K'')$ is nef and big on $M^{\sharp}$ for any sufficiently ample line bundle $H$ on $W$, so
$H^q(M^{\sharp},H+E+tA)=0$ for $q>0$, $t\ge\nomb -1$.
Likewise $H+E''+tA-D^- -(K^{\sharp}+\frac12E'')$ is nef and big on $M^{\sharp}$ and $H^q(M^{\sharp},H+E''+tA-D^-)=0$ for $q>0$, $t\ge -1$.

Using the exact sequences $0 \to \SO(H+tA) \to \SO(H+E+tA) \to \SO_E(-2) \to 0$ and $0 \to \SO(H+tA-D^-) \to \SO(H+E''+tA-D^-) \to \SO_{E''}(-2) \to\nomb 0$,
we infer $H^q(M^{\sharp},H+tA)=0$ and $H^q(M^{\sharp},H+tA-D^-)=0$ for $q>0$, $t\ge\nomb -1$.
Hence $H^q(D^-,tA)=0$ for $q>0$, $t\ge -1$, which implies
$\chi(D^-,tA)=0$ for $t=-1$ and ${}=1$ for $t=0$.
Therefore $\chi(D^-,tA)=(t+1)(dt+2)/2$ for $d=A^2D$.

For any line bundle $F$ on $M^{\sharp}$ such that $F_{E_i}=\SO(j)$ for each component $E_i$ of $E$, we have
$\chi(D^-,F-E)-\chi(D^-,F)=\chi(E,F-D^-)-\chi(E,F)=\sum'_i(\chi(E_i,\SO(j-1))-\chi(E_i,\SO(j)))=-\mu'(j+1)$, where $\mu'$ is the number of components of $E'$.
Setting $B=L^{\sharp}-2A\in\Pic(M^{\sharp})$, we get
$\chi(D^-,2tB)=\chi(D^-,tA)-\sum_{j=0}^{t-1}\mu'(2j+1)=(t+1)(dt+2)/2-\mu't^2$
by applying this formula successively, since $A=2B+E$ in $\Pic(D^-)$.
This is true for all half integers, so
$\chi(D^-,-B)=\frac12-\frac d8-\frac{\mu'}4$ and
$\chi(D^-,-3B)=\frac38{d}-\frac12-\frac94\mu'$.

On the other hand, $A-D$ is $\rho$-ample and $H-B-D^- -(K^{\sharp}+\frac12E')=\pi^*(H+A-D)$ and $H-B-K^{\sharp}=\pi^*(H+A)$ are nef and big on $M^{\sharp}$, so we get
$h^q(M^{\sharp},H-B-D^-)=h^q(M^{\sharp},H-B)=0$ for $q>0$, hence $\chi(D^-,-B)=0$ for $q>0$.
Likewise $H-B-A-K^{\sharp}$ and $H-B-A-D^- -(K^{\sharp}+\frac12E')$ are nef and big on $M^{\sharp}$, so
$0=\chi(D^-,-B-A)=\chi(D^-,-3B-E)$, hence
$\chi(D^-,-3B)=\chi(E,-3B)-\chi(E,-3B-D^-)=\sum'_i(\chi(E_i,\SO(-3))-\chi(E_i,\SO(-4)))=-2\mu'$.
Now, combining with the above formula, we get $d+2\mu'=4$ and $3d=4+2\mu'$, hence $d=2$ and $\mu'=1$.

As in [{\bf F4};(4.6.4;f)], we infer $h^0(D,A)=h^0(D^-,A)=\chi(D^-,A)=d+\nomb 2=4$ and $\Delta(D,A)=0$, so $(D, A)$ is a hyperquadric.
As in (3.6.3;f), $M''$ has no singularity of the type (2.3;3) on the smooth locus of $D$, so $\mu'=1$ implies that $D$ is a cone over a quadric curve and $\pi(E')$ is the vertex $v$ of it.
$D^{\sharp}$ is the blow-up of $D$ at $v$, so $D^{\sharp}\simeq\Sigma_2$ and $E\cap D^{\sharp}$ is the $(-2)$-curve on it.
Since $K^{\sharp}=K''+\frac12 E$ is numerically equivalent to $-\frac32 A+\frac12 E$ on $D^{\sharp}$,
the normal bundle $[D^{\sharp}]_{D^{\sharp}}$ of $D^{\sharp}$ is $-A+Y$ by the adjunction formula,
where $Y$ is a fiber of $\Sigma_2\to \BP^1$.
Hence $D^{\sharp}$ can be blown down, and then the image of $E$ can be blown down succesively to a smooth point.
The result is nothing but $W$, and we conclude that $\rho(D)$ is a smooth point on $W$, as in [{\bf F4};(4.6.4;h)].

(3.8.$\infty$) {\it Summary of the case $\tau''=\frac32$.}

An extremal ray $R$ such that $(K''+\frac32A)R=0$ is of the type (3.8.0.0), (3.8.0.1.0), (3.8.0.2), (3.8.0.3), (3.8.1) or (3.8.3).
If there is a ray which is of one of the five fibration types, then the structure of $(M'', A)$ is almost determined and $\ke(M, L)=-3/5$ in these cases.

Suppose that every such ray is of the type (3.8.3).
As in the cases of first and second reductions, the exceptional divisors of different rays are disjoint with each other, so we can blow down them simultaneously to smooth points, $\beta: M'' \lra M^\flat$.
By construction $M^\flat$ has no worse singularities than $M''$ as in (3.6.3), and we can get rid of this case by replacing $(M'',A)$ by $(M^\flat,A^\flat)$.

(3.9) The case $\tau''=4/3$.

We have $(3K''+4A)R=0$ for some ray $R$.
Set $U=5A-2L''\in\Pic(M'')$.
Then $3UR=AR$ and $(K''+4U)R=0$.
This implies $(M'', U)\simeq(\BP^3,\SO(1))$ by (1.6).
Moreover $A=\SO(3)$, $L''=\SO(7)$ and $\ke(M, L)=-4/7$.
Clearly there are many examples of this type.

(3.10) The case $\tau''=5/4$.

We have $(4K''+5A)R=0$ for some ray $R$.
Let $\rho: M'' \to W$ be the contraction morphism of $R$.
Set $U=2L''-4A\in\Pic(M'')$.
Then $AR=2UR$, $K''R=-\frac52 UR$ and $L''R=\frac92 UR$.

(3.10.0) The case $\dm W=0$.

$(M'',U)$ is a polarized variety with nef value $\tau=5/2$.
Set $B:=L^{\sharp}-2A\in\Pic(M^{\sharp})$.
Then $2B+E=U$ in $\Pic(M^{\sharp})$, so we infer that $(M'',U)$ is the projective cone over the Veronese surface $(\BP^2,\SO(2))$ by the same argument in (3.5.0).
Clearly this case occurs really.

(3.10.i) The cases $1\le\dm W<3$.

Any general fiber $F$ of $\rho$ is smooth.
These cases are ruled out as in (3.5).

(3.10.3) The case $\dm W=3$. $\rho$ is birational.

This case is ruled out as in (3.5.3).
We replace $A$ by $U$ here.

(3.11) The case $\tau''=7/6$.

We have $(6K''+7A)R=0$ for some ray $R$.
Set $U=2L''-4A\in\Pic(M'')$.
Then $AR=3UR$ and $K''R=-\frac72 UR$.
But this contradicts (1.6), thus this case is ruled out.

(3.12) The next possible value for $\tau''$ is $\tau''=1$, but this corresponds to cases $\ke\ge -1/2$.
Thus, the above arguments are enough to classify the cases $\ke<-1/2$.
Here we present some partial results in case $\tau''=1$.

By the Fibration Theorem, there is a morphism $g: M'' \to W$ onto a normal variety $W$ such that $g_*\SO_{M''}=\SO_W$ and $K''+A=g^*H$ for some ample $\BQ$-bundle $H$ on $W$.

(3.12.0) The case $\dm W=0$.

By duality we have $h^3(M'',-A)=h^0(M'',\omega\otimes A)\le 1$,
where $\omega$ is the canonical sheaf of $M''$.

If $h^3(M'',-A)=1$, then $K''$ is invertible and $K''=-A$.
Thus $M''$ has no singulariy of the type (2.3; 3) and $M''$ itself is a Fano threefold (possibly with some singularities of the type (2.3; 2)).

But it actually occurs that $h^3(M'',-A)=0$,
if $M''$ has a singularity of the type (2.3; 3) and $\omega$ is not invertible.
We have $\chi(M'',-A)=0$ and $\chi(M'',\SO)=1$, so
by the Riemann-Roch Theorem we can set
$\chi(M'',tA)=\frac d6 t^3+at^2+bt+1$ for some numbers $a$, $b$ and $d=A^3$.
Moreover $a=d/4$ since $K''$ is numerically $-A$.
Combining these observations we get
$\chi(M'',tA)=(\frac d{12}(2t+1)t+1)(t+1)$.

Let $\pi: M^{\sharp} \to M''$ be the factorial stage and let $\mu$ be the number of components of the exceptional divisor $E$ of $\pi$.
Then the canonical bundle $K^{\sharp}$ of $M^{\sharp}$ is $K''+\frac12E=\frac12 E-A$, so $E=2(K^{\sharp}+A)$ and
$\chi(M^{\sharp},tE)=\chi(M^{\sharp},\SO)+\sum_{j=1}^t\chi(E,\SO(-2j))=1+\mu\sum_j(2j-1)(2j-2)/2=1+\mu t(t-1)(4t+1)/6$.
Setting $B:=K^{\sharp}+A$, we have $H^q(M^{\sharp},B)=0$ for $q>0$ by the Vanishing Theorem, while $h^0(M^{\sharp},B)=h^3(M^{\sharp},-A)=h^3(M'',-A)=0$, hence
$0=\chi(M^{\sharp},B)=\chi(M^{\sharp},\frac12 E)=1-\frac\mu 8$, so $\mu=8$.
Thus $\chi(M^{\sharp},sB)=(s+1)(s-1)(2s-3)/3$.

Any way, since $E\in\vert 2B\vert$, there is a double covering $f: N^{\sharp} \to M^{\sharp}$ branched along $E$ such that $f_*\SO_{N^{\sharp}}=\SO_{M^{\sharp}}\oplus\SO(-B)$.
Note that $f^*E=2D$ for some divisor $D$ on $N^{\sharp}$ such that $D\simeq E$.
Moreover $[D_j]_{D_j}=\SO(-1)$ for each component $D_j\simeq\BP^2$ of $D$.
Hence $D$ can be blown down to smooth points; $N^{\sharp} \to N$.
The sheet changing involution of the double cover $f$ induces an involution $\theta$ of $N$ having exactly eight fixed points, and the natural map $N \to M''$ is the quotient map.
The canonical bundle of $N^{\sharp}$ is $f^*(K^{\sharp}+B)=f^*(2B-A)=2D-f^*A$ and hence the canonical bundle of $N$ is $-A_N$.
Thus $N$ is a Fano threefold.

Sano [{\bf S}] classifies possible types of this Fano threefold $N$, and gives various interesting examples of such pairs $(N,\theta)$.
There he says that one can easily classify pairs with ``a very tiresome work'', details of which are omitted.
We should also remark that some of such pairs cannot appear as the second reduction of a polarized manifold with $\ke=-1/2$.

(3.12.1) The case $\dm W=1$.

Any general fiber $F$ of $g$ is smooth and $K_F+A_F=0$.
Thus $(F,A_F)$ is a Del Pezzo surface.
The types of possible singular fibers are not yet classified;
unlike [{\bf F3}], $M''$ may have singularities.

(3.12.2) The case $\dm W=2$.

Any general fiber $F$ of $g$ is $\BP^1$, and $AF=2$, $L''F=4$.

(3.12.3) The case $\dm W=3$.

In this case $K''+A$ is nef big and $\ke(M,L)>-1/2$.

(3.$\infty$) {\it Summary.}
For $t=-\ke(M,L)$, the Iitaka dimension of the $\BQ$-bundle $K+tL$ will be called {\it the adjoint Kodaira dimension} of $(M,L)$, denoted by $\alpha\kappa(M,L)$.
Then
smooth polarized threefolds with $\ke(M,L)<-1/2$ can be classified as follows:
\dnl
$\ke=-4$; $\alpha\kappa=0$. $(M,L)\simeq(\BP^3,\SO(1))$.
\dnl
$\ke=-3$; $\alpha\kappa=0$. $(M,L)$ is a hyperquadric $(\BQ^3,\SO(1))$ in $\BP^4$.
\nl\indent\qquad $\alpha\kappa=1$. $(M,L)$ is a scroll over a curve.
\dnl
$\ke=-2$; $\alpha\kappa=0$. $(M,L)$ is a Del Pezzo threefold.
\nl\indent\qquad $\alpha\kappa=1$. $(M,L)$ is a hyperquadric fibration over a curve.
\nl\indent\qquad $\alpha\kappa=2$. $(M,L)$ is a scroll over a surface.
\dnl
$\ke=-\frac32$; $\alpha\kappa=0$. The first reduction $(M',L')$ is $(\BQ^3,\SO(2))$.
\nl\indent\qquad $\alpha\kappa=1$. $(M',L')$ is a Veronese fibration over a curve, i.e., any \nl
\indent\qquad\qquad\qquad general fiber  is $(\BP^2,\SO(2))$.
\dnl
$\ke=-\frac43$; $\alpha\kappa=0$. $(M',L')\simeq(\BP^3,\SO(3))$.
\dnl
$\ke=-1$; $\alpha\kappa=0$. $(M',L')$ is a Fano threefold with canonical bundle $-L'$.
\nl\indent\qquad $\alpha\kappa=1$. $(M',L')$ is a Del Pezzo fibration over a curve.
\nl\indent\qquad $\alpha\kappa=2$. $(M',L')$ is a hyperquadric fibration over a surface.
\dnl
$\ke=-\frac45$; $\alpha\kappa=0$. The second reduction $(M'',A)$ is $(\BP^3,\SO(1)$, and \nl
\indent\qquad\qquad\qquad $L''=\SO(5)$, cf. (3.2).
\dnl
$\ke=-\frac34$; $\alpha\kappa=0$. $(M'',A)$ is a hyperquadric.
\nl\indent\qquad $\alpha\kappa=1$. $(M'',A)$ is a scroll over a curve.
\dnl
$\ke=-\frac57$; $\alpha\kappa=0$. $(M'',A)$ is the projective cone over $(\BP^2,\SO(2))$, cf. (3.5).
\dnl
From now on, one may need some contractions of the type (3.6.3) to obtain $(M'',A)$.
\dnl
$\ke=-\frac23$; $\alpha\kappa=0$. $(M'',A)$ is a Del Pezzo threefold, cf. (3.6.0).
\nl\indent\qquad $\alpha\kappa=1$. $(M'',A)$ is a hyperquadric fibration over a curve, cf. (3.6.1).
\nl\indent\qquad $\alpha\kappa=2$. $(M'',A)$ is a scroll over a surface, cf. (3.6.2).
\dnl
From now on, one may need further some contractions of the type (3.8.3) to obtain $(M'',A)$.
\dnl
$\ke=-\frac35$; $\alpha\kappa=0$. $(M'',A)$ is one of the types described in (3.8.0), (3.8.0.1.0), \nl
\indent\qquad\qquad\qquad (3.8.0.2) or (3.8.0.3).
\nl\indent\qquad $\alpha\kappa=1$. $(M'',A)$ is a Veronese fibration, cf. (3.8.1).
\dnl
$\ke=-\frac47$; $\alpha\kappa=0$. $(M'', A)$ is $(\BP^3,\SO(3))$.
\dnl
$\ke=-\frac59$; $\alpha\kappa=0$. $(M'',A)$ is of the form $(M'',2U)$, where $(M'',U)$ is the \nl
\indent\qquad\qquad\qquad projective cone over $(\BP^2,\SO(2))$, cf. (3.10).

\enddocument